\documentclass[12pt,a4paper]{article}
\usepackage{hyperref}
\usepackage{amsfonts}
\usepackage{mathrsfs}
\usepackage{comment}
\usepackage{amsmath}
\usepackage{amssymb}
\usepackage{amsthm}
\usepackage{amscd}
\usepackage{graphicx}
\usepackage{indentfirst}
\usepackage[all]{xy}
\usepackage{titlesec}
\usepackage[top=20mm, bottom=20mm, left=20mm, right=20mm]{geometry}
\usepackage{enumerate}
\usepackage{bm}
\usepackage{enumitem}
\usepackage{color}
\usepackage{bbm}
\newtheorem{theorem}{Theorem}[]
\newtheorem{lemma}[theorem]{Lemma}

\newtheorem{conjecture}[theorem]{Conjecture}

\newtheorem{remark}[theorem]{Remark}

\newcommand{\ma}{\mathcal}

\newcommand{\s}{\subseteq}

\newcommand{\fr}{\frac}

\newcommand{\mat}{\mathbbm}

\begin{document}

\title{Maximum subsets of $\mathbb{F}^n_q$ containing no right angles}
\date{}
\author{
Gennian Ge$^{\text{a,}}$\thanks{Email address: gnge@zju.edu.cn. Research supported by the National Natural Science Foundation of China under Grant Nos. 11431003, 61571310 and 11971325, Beijing Scholars Program, Beijing Hundreds of Leading Talents Training Project of Science and Technology, and Beijing Municipal Natural Science Foundation.},
Chong Shangguan$^{\text{a,b}}$\thanks{Corresponding author. Email address: theoreming@163.com.}\\
\footnotesize $^{\text{a}}$ School of Mathematical Sciences, Capital Normal University, Beijing, 100048, China.\\
\footnotesize $^{\text{b}}$ School of Mathematical Sciences, Zhejiang University, Hangzhou 310027, China.\\
}

\maketitle

\begin{abstract}
    Recently, Croot, Lev, and Pach (Ann. of Math., 185:331--337, 2017.) and Ellenberg and Gijswijt (Ann. of Math., 185:339--443, 2017.) developed a new polynomial method and used it to prove upper bounds for three-term arithmetic progression free sets in $\mathbb{Z}_4^n$ and $\mathbb{F}_3^n$, respectively.
    Their approach was later summarized by Tao and is now known as the slice rank method.
    In this paper, we apply this method to obtain a new upper bound on the cardinality of subsets of $\mathbb{F}^n_q$ which contain no right angles.

    More precisely, let $q$ be a fixed odd prime power and $x\cdot y$ be the standard inner product of two vectors $x,y\in\mathbb{F}_q^n$, we prove that the maximum cardinality of a subset $A\s\mathbb{F}_q^n$ without three distinct elements $x,y,z\in A$ satisfying $(z-x)\cdot (y-x)=0$ is at most $\binom{n+q}{q-1}+3$.
    For sufficiently large $n$, our result significantly improves the previous upper bound of Bennett (European J.
    Combin., 70:155--163, 2018.), who showed that $|A|=\mathcal{O}(q^{\frac{n+2}{3}})$.
\end{abstract}



\section{Introduction}
\noindent The recent breakthrough results of Croot, Lev, and Pach \cite{croot2016progression} and Ellenberg and Gijswijt \cite{ellenberg2016large} showed respectively that three-term arithmetic progression free sets in $\mathbb{Z}_4^n$ and $\mathbb{F}_3^n$ are exponentially small.
The main idea in their proofs was an application of a novel polynomial method, which was later summarized by Tao \cite{tao} as a principle which compares the slice ranks (see Section 2 below) of different multivariate functions.

The problems investigated in \cite{croot2016progression,ellenberg2016large} can be viewed as extremal problems defined over the finite fields which forbid the existence of certain configurations.
In order to deal with this type of problems by the polynomial method, the essential idea is to characterize the property of having no such configurations by an equation between two carefully chosen functions.
For example, let the forbidden configuration be the three-term arithmetic progression, then a subset $A\s\mathbb{F}_3^n$ contains no such configurations if and only if for three elements $x,y,z\in A$, $x+y+z=0^{n}$ if and only if $x=y=z$.
In \cite{tao} this observation is described as an equation between two $\mathbb{F}_3$-valued functions defined on $A\times A\times A$:

\begin{equation}\label{principal}
    \begin{aligned}
\mathbbm{1}_{x+y+z=0^n}=\mathbbm{1}_{x=y=z},\\
    \end{aligned}
\end{equation}

\noindent where for any $I$ (here $I$ can be an equation, an event, etc.), $\mat{1}_I$ is the binary indicator function such that $\mat{1}_I=1$ if $I$ holds, otherwise $\mathbbm{1}_I=0$.
Thus by comparing the slice ranks of the two functions on the left and right hand sides of (\ref{principal}), one can possibly bound the cardinality of $A$ from above.

Since the works of \cite{croot2016progression,ellenberg2016large,tao}, the slice rank method has enjoyed immense success and it has been applied to many different problems, for example, the tri-colored sum-free sets \cite{Alonsomething,tricolored,Kleinberg2016The} and the sunflower-free sets \cite{Alonsomething,Naslund2016Upper}.
The interested reader is referred to \cite{grochow2018new} for a survey on more applications and to \cite{tricolored,ChriSTOC,tao2} for further developments of this method.

The goal of this paper is to apply the slice rank method to deal with the following extremal problem.
%defined on the finite fields.
Let $q$ be a prime power and $V:=\mathbb{F}_q^n$ be an $n$-dimensional vector space on the finite field $\mathbb{F}_q$.
For any two vectors $x=(x_1,\ldots,x_n),y=(y_1,\ldots,y_n)\in V$, let $x\cdot y$ be the standard inner product such that $x\cdot y=\sum_{i=1}^nx_iy_i\in\mathbb{F}_q$.
A right angle of $V$ is a subset of three {\it distinct} elements $x,y,z\in V$ such that $(z-x)\cdot(y-x)=0$.
We are interested in the maximum cardinality of subsets of $V$, in which no right angles are contained.
This problem is a finite field version of the Erd{\H{o}}s-Falconer type problem, which was originally defined in the setting of Euclidean spaces \cite{zhuangbi2,Falconer1992On,Euclidean} and asked for the smallest $d$ for which any compact set in $\mathbb{R}^n$ with Hausdorff dimension larger than $d$ contains three points forming an angle $\alpha$, for given $n$ and $\alpha$.
The reader is referred to \cite{Erdosfalconer} for more finite field analogs of problems in Euclidean space.

Let $R(n,q)$ denote the maximum cardinality of subsets of $\mathbb{F}_q^n$ which contain no right angles.
In \cite{bennett2015right} Bennett showed that for any $q$ and $n$

\begin{equation}\label{1}
  \begin{aligned}
    R(n,q)=\ma{O}(q^{\fr{n+2}{3}}).
  \end{aligned}
\end{equation}

\noindent For fixed $n$ and large $q$, the above bound was improved slightly by Pham, Sang and Tardos \cite{phamright} to

\begin{equation}\label{2}
  \begin{aligned}
&R(n,q)=\ma{O}(nq^{\fr{n+1}{3}}\log q), &n=3k-1,~k\ge 1,~\text{and}\\
&R(n,q)=c_nq^{\fr{n+1}{3}+\fr{k}{3k-1}}(\log q)^{\fr{k}{3k-1}},&n=3k~\text{or}~n=3k+1,~k\ge 1,\\
  \end{aligned}
\end{equation}

\noindent where $c_n$ is a positive constant depending on $n$.

It is easy to verify that if a subset $A\s\mathbb{F}_q^n$ contains no right angles, then for any three (not necessarily distinct) elements $x,y,z$ of $A$, $(z-x)\cdot(y-x)=0$ if and only if at least one of the following holds: $(i)$ $z=x$ or $y=x$; $(ii)$ $z=y$ and $(y-x)\cdot (y-x)=0$.
The property of containing no right angles can also be characterized by an equation between two $\mathbb{F}_q$-valued polynomials, as in (\ref{principal}).
Thus by applying the slice rank method one can obtain a new upper bound on $R(n,q)$, as stated below.

\begin{theorem}\label{rightang}
    Let $q$ be a fixed odd prime power and $A\s\mathbb{F}_q^n$ be a subset containing no three distinct elements $x,y,z\in A$ such that $(z-x)\cdot(y-x)=0$. Then $|A|\le \binom{n+q}{q-1}+3$.
\end{theorem}

Clearly, for fixed odd $q$ and $n\rightarrow\infty$, (\ref{1}) and (\ref{2}) are both exponential functions of $n$, whereas the new bound in Theorem \ref{rightang} is only a polynomial function of $n$.
Thus the previous results are significantly improved.
For example, for $q=3$, $R(n,3)$ is improved from $\ma{O}(3^{\fr{n+2}{3}})$ to $\binom{n+3}{2}+3$.
In this sense Theorem \ref{rightang} is of interest since all previous applications of the slice rank method resulted in exponential bounds (e.g., \cite{croot2016progression,ellenberg2016large,Kleinberg2016The,Naslund2016Upper}).
Note that since the standard orthonormal basis $\{e_1,\ldots,e_n\}$ contains no right angles, for any prime power $q$ we have that $R(n,q)\ge n$.
In general, we believe that the following conjecture is true.

\begin{conjecture}\label{conjecture}
For any fixed prime power $q$, $R(n,q)=\Theta(n^{q-1})$ as $n\rightarrow\infty$.
\end{conjecture}

The rest of this paper is organised as follows. In Section 2 we review the necessary terminologies introduced in \cite{tao} and in Section 3 we present the proof of Theorem \ref{rightang}.

\section{The slice rank and a lemma of Tao}

\noindent Let $\mathbb{F}$ be an arbitrary field and $A$ be a finite set.
Tao \cite{tao} in his blog introduced a rank notion (which was later termed the slice rank) for $A^k\rightarrow\mathbb{F}$ functions with $k$ variables $x_1,\ldots,x_k$, as defined below.

\begin{itemize}
  \item Constant functions have slice rank zero.

  \item A non-constant function $h:A^k\rightarrow\mathbb{F}$ has slice rank one if and only if it can be expressed as $$h(x_1,\ldots,x_k)=f(x_i)g(x_1,\ldots,x_{i-1},x_{i+1},\ldots,x_k)$$ for some integer $1\le i\le k$ and functions $f:A\rightarrow\mathbb{F}$ and $g:A^{k-1}\rightarrow\mathbb{F}$.

  \item In general, for any function $T:A^k\rightarrow\mathbb{F}$, its {\it slice rank} ${\rm{sr}}(T)$ is the smallest nonnegative integer $r$ such that $T$ can be written as the sum of $r$ functions of slice rank one.
\end{itemize}

For example, if $k=2$, then a slice rank one function has the form $T(x,y)=f(x)g(y)$ for some $f,g:A\rightarrow\mathbb{F}$;
if $k=3$, then the slice rank one functions have one of the following forms:
$$T_1(x,y,z)=f_1(x)g_1(y,z)\text{ or }T_2(x,y,z)=f_2(y)g_2(x,z)\text{ or }T_3(x,y,z)=f_3(z)g_3(x,y)$$
for some $f_i:A\rightarrow\mathbb{F}$ and $g_i:A^2\rightarrow\mathbb{F}$, $1\le i\le 3$.
It can be easily seen that the linear combination of $r$ functions of slice rank one produces a function of slice rank at most $r$.
Consequently ${\rm{sr}}(\cdot)$ is subadditive, i.e., for $A^k\rightarrow \mathbb{F}$ functions $T_1,\ldots,T_r$,
\begin{equation}\label{subadditivity}
  \begin{aligned}
    {\rm{sr}}(\sum_{i=1}^r T_i)\le \sum_{i=1}^r {\rm{sr}}(T_i).
  \end{aligned}
\end{equation}

To prove Theorem \ref{rightang} we will make use of the following lemma, which is a special case of a result in \cite{tao}.

\begin{lemma}[see Lemma 1 of \cite{tao}]\label{tao}
   Let $A$ be a finite set.
   The the slice rank of the function
   \begin{equation*}
     \begin{aligned}
       (x,y,z)\in A\times A\times A\mapsto \mathbbm{1}_{x=y=z}
     \end{aligned}
   \end{equation*}

   \noindent is $|A|$ over any field.
\end{lemma}

\section{Large subsets of $\mathbb{F}_q^n$ containing no right angles}

\begin{proof}[\textbf{Proof of Theorem \ref{rightang}}] Define a function $T:A\times A\times A\rightarrow\mathbb{F}_q$ as follows

\begin{equation*}\label{rightangleft}
    \begin{aligned}
           T(x,y,z)=\mathbbm{1}_{y\neq z}\cdot\big((y-x)\cdot (z-x)\big)^{q-1}.
    \end{aligned}
\end{equation*}

\noindent Observe that since $A$ contains no right angles, then $(y-x)\cdot (z-x)\neq0$ for distinct $x,y,z\in A$, which implies that $\big((y-x)\cdot (z-x)\big)^{q-1}=1$ and hence

\begin{equation*}\label{3}T(x,y,z)=
        \begin{cases}
        0, & x=y\text{ or }x=z\text{ or }y=z,\\
        1, & otherwise.
        \end{cases}
    \end{equation*}

\noindent It is easy to verify that
$$T(x,y,z)=\mathbbm{1}_{(x\neq y)\wedge (x\neq z)\wedge (y\neq z)}=2\cdot\mathbbm{1}_{x=y=z}-\mathbbm{1}_{x=y}-\mathbbm{1}_{x=z}-\mathbbm{1}_{y=z}+1.$$

\noindent Since $2\nmid q$, one can write $\mathbbm{1}_{x=y=z}=\fr{1}{2}\big(T(x,y,z)+\mathbbm{1}_{x=y}+\mathbbm{1}_{x=z}+\mathbbm{1}_{y=z}-1\big)$.
By Lemma \ref{tao} and (\ref{subadditivity}) it holds that
%$|A|={\rm{sr}}(\mathbbm{1}_{x=y=z})\le {\rm{sr}}(T) +3$ and hence

\begin{equation}\label{lowerbound}
  \begin{aligned}
        |A|={\rm{sr}}(\mathbbm{1}_{x=y=z})\le {\rm{sr}}(T) +3.
        %{\rm{sr}}(T)\ge |A|-3,
        %={\rm{sr}}(2\cdot\mathbbm{1}_{x=y=z}-\mathbbm{1}_{x=y}-\mathbbm{1}_{x=z}-\mathbbm{1}_{y=z}+1)\ge |A|-3,
  \end{aligned}
\end{equation}

\noindent On the other hand, notice that %by expanding $\big((y-x)\cdot (z-x)\big)^{q-1}$ as the sum of monomials

\begin{equation*}
  \begin{aligned}
    \big((y-x)\cdot (z-x)\big)^{q-1}&=\big(\sum_{i=1}^ny_iz_i+\sum_{i=1}^nx_i^2-x_1(y_1+z_1)-\cdots-x_n(y_n+z_n)\big)^{q-1}\\
    &=\sum_{l+m+\sum_{i=1}^n k_i=q-1}a_{l,m,k_1,\ldots,k_n}
    \big(\sum_{i=1}^ny_iz_i\big)^{l}\big(\sum_{i=1}^nx_i^2\big)^{m}\prod_{i=1}^n\big(x_i(y_i+z_i)\big)^{k_i},\\
    %\big(x_1(y_1+z_1)\big)^{k_1}\cdots\big(x_n(y_n+z_n)\big)^{k_n},\\
  \end{aligned}
\end{equation*}

\noindent where $l,m,k_1,\cdots,k_n$ are nonnegative integers summing up to $q-1$ and $a_{l,m,k_1,\ldots,k_n}\in\mathbb{F}_q$ is the coefficient of the corresponding monomial. Thus

\begin{equation*}
    \begin{aligned}
        T(x,y,z)&=\mathbbm{1}_{y\neq z}\cdot
        \sum_{l+m+\sum_{i=1}^n k_i=q-1}a_{l,m,k_1,\ldots,k_n}
    \big(\sum_{i=1}^ny_iz_i\big)^{l}\big(\sum_{i=1}^nx_i^2\big)^{m}\prod_{i=1}^n\big(x_i(y_i+z_i)\big)^{k_i}\\
    &=\sum_{l+m+\sum_{i=1}^n k_i=q-1}a_{l,m,k_1,\ldots,k_n}
    \Big(\big(\sum_{i=1}^nx_i^2\big)^{m}\prod_{i=1}^nx_i^{k_i}\Big)\cdot
    \Big(\mathbbm{1}_{y\neq z}\cdot\big(\sum_{i=1}^ny_iz_i\big)^{l}\prod_{i=1}^n(y_i+z_i)^{k_i} \Big)\\
    &:=\sum_{l+m+\sum_{i=1}^n k_i=q-1}a_{l,m,k_1,\ldots,k_n} f_{m,k_1,\ldots,k_n}(x)g_{l,k_1,\ldots,k_n}(y,z),\\
    \end{aligned}
\end{equation*}

\noindent where we set
$f_{m,k_1,\ldots,k_n}(x)=\big(\sum_{i=1}^nx_i^2\big)^{m}\prod_{i=1}^nx_i^{k_i}\text{~and~}g_{l,k_1,\ldots,k_n}(y,z)=\mathbbm{1}_{y\neq z}\cdot\big(\sum_{i=1}^ny_iz_i\big)^{l}\prod_{i=1}^n(y_i+z_i)^{k_i}.$
Observe that ${\rm{sr}}\big(a_{l,m,k_1,\ldots,k_n}f_{m,k_1,\ldots,k_n}(x)g_{l,k_1,\ldots,k_n}(y,z)\big)\le 1$ (possibly zero), then

\begin{equation}\label{upperbound}
  \begin{aligned}
    {\rm{sr}}(T)\le\binom{n+q}{q-1},
  \end{aligned}
\end{equation}

\noindent since there are $\binom{n+q}{q-1}$ different choices of nonnegative integers $l,m,k_1,\ldots,k_n$ such that $l+m+\sum_{i=1}^n k_i=q-1$ and $T$ can be expressed as the sum of at most $\binom{n+q}{q-1}$ slice rank one functions.
Combining (\ref{lowerbound}) and (\ref{upperbound}) it follows that

\begin{equation*}
  |A|-3\le {\rm{sr}}(T) \le \binom{n+q}{q-1},
\end{equation*}

\noindent which proves the theorem.
\end{proof}

\begin{remark}
  After the paper was submitted, the authors noticed that the result of Theorem \ref{rightang} was further generalized by Naslund \cite{naslund2017partition} to the $k$-right-corner-free sets, where a $k$-right-corner is a subset of $k+1$ distinct vectors $x_1,\ldots,x_{k+1}\in\mathbb{F}_q^n$ such that $x_1-x_{k+1},\ldots,x_k-x_{k+1}$ are pairwise orthogonal.
  In \cite{naslund2017partition} it was shown that for any prime power $q$ with characteristic $p$ and $p>k$, it holds that any $A\s\mathbb{F}_q^n$ with $|A|>(k+1)\cdot\binom{n+(k-1)q}{(k-1)(q-1)}$ contains a $k$-right-corner.
  More interestingly, the proof of this result used a more general version of the slice rank, which was termed the {\it partition rank} \cite{naslund2017partition}.
\end{remark}

\vspace{10pt}

\noindent\textbf{Acknowledgement.} The authors are grateful to the three anonymous referees for their careful reading and constructive comments which are very helpful to the improvement of the presentation of this paper, and they would also like to thank Xiangliang Kong, Jingxue Ma and Yiwei Zhang for helpful discussions.
\bibliographystyle{plain}
\bibliography{rightang}

\end{document}